\documentclass[10pt]{article}

\usepackage{url}

\usepackage{authblk}
\usepackage{mathtools}
\usepackage{amsfonts}
\usepackage{graphicx}
\usepackage{verbatim}
\usepackage{mathrsfs}
\usepackage{subfigure}
\usepackage{fancyhdr}
\usepackage{latexsym}
\usepackage{graphicx}
\usepackage{dsfont}
\usepackage{fancybox}
\usepackage{nicefrac}
\usepackage{framed, color}
\usepackage{floatrow}
\definecolor{shadecolor}{rgb}{1, 0.8, 0.3}
\usepackage {blindtext}
\usepackage{geometry}
\usepackage{hyperref}
\usepackage{enumitem}

\normalfont
\usepackage[T1]{fontenc}
\usepackage{relsize}
\usepackage[utf8]{inputenc}
\usepackage[english]{babel}
\usepackage{enumitem,graphicx,appendix}

\usepackage{xargs}
\usepackage{color}

\def\AA{\mathcal{A}}
\def\BB{{\mathcal B}}
\def\FF{{\mathcal F}}

\def\II{{\mathcal I}}

\def\NN{{\mathcal N}}
\def\XX{{\mathcal X}}
\def\YY{{\mathcal Y}}

\def\complex{\mathbb{C}}

\def\MM{\mathcal{G}}

\def\proj{{\mathbb P}}
\def\Id{{\mathbb I}}
\def\half{{\frac12}}

\def\defeq{\displaystyle\mathrel{\mathop=^{\scriptscriptstyle\rm def}}}

\def\bdot{\hbox{\bf .}}

\def\id{{\mathbb I}}
\def\eps{{\varepsilon}}

\def\chebb{{Chebyshev}}
\def\cheb{{Chebyshev\ }}

\newtheorem{theorem}{Theorem}

\newtheorem{lemma}{Lemma}
\newtheorem{proposition}{Proposition}

\usepackage{hyperref}
\hypersetup{
    colorlinks=true,
    linkcolor=blue,
    filecolor=magenta,
    urlcolor=cyan,
}

\newtheorem{definition}{Definition}

\SetLabelAlign{Center}{\hfil \mbox({#1}) \hfil}
\newfloatcommand{capbtabbox}{table}[][\FBwidth]

\makeatletter
\newcommand{\tpitchfork}{%
  \vbox{
    \baselineskip\z@skip
    \lineskip-.52ex
    \lineskiplimit\maxdimen
    \m@th
    \ialign{##\crcr\hidewidth\smash{$-$}\hidewidth\crcr$\pitchfork$\crcr}
  }%
}
\makeatother

\title{
Branches and bifurcations of  ejection-collision orbits \\
in the planar circular restricted three body problem
}

\author[1]{Gianni Arioli \thanks{
G.A. partially supported by PRIN project 2022 "Partial differential equations and related geometric-functional inequalities", financially supported by the EU, in the framework of the "Next Generation EU initiative".
Email: {\tt gianni.arioli@polimi.it}}}

\author[2]{J.D. Mireles James \thanks{J.M.J partially supported by NSF 
grant DMS - 2307987.
Email: {\tt jmirelesjames@fau.edu}}}

\affil[1]{Politecnico di Milano, Department of Mathematics}

\affil[2]{Florida Atlantic University, Department of Mathematical Sciences}
\date{\today}

\begin{document}

\maketitle

\begin{abstract}
The goal of this paper it to
prove existence theorems for
one parameter families (branches) of ejection-collision orbits 
 in the planar circular restricted three body problem (CRTBP),
 and to study some of bifurcations of these branches.
The CRTBP considers the dynamics of an infinitesimal particle
moving in the gravitational field of two massive primary bodies.
The motion of the primaries assumed to be circular, and we 
study ejection-collision orbits where the infinitesimal body is 
ejected from one primary and collides with the other
(as opposed to more local ejections-collisions
where the infinitesimal body collides with a single primary body 
in both forward and backward time).  We consider branches of 
ejection-collision orbits which are 
(i) parameterized by the Jacobi integral (energy like quantity conserved 
by the CRTBP) with the masses of the primaries fixed, and 
(ii) parameterized by the mass ratio of the primary bodies with energy fixed.
The method of proof is constructive and computer assisted, hence can be
applied in non-perturbative settings and (potentially) to other conservative
systems of differential equations. The main requirement  
is that the system should admit a change of coordinates which 
regularizes the collision singularities.  In the planar CRTBP,
the necessary regularization is  provided by the classical 
Levi-Civita transformation. 
\end{abstract}

\section{Introduction} \label{sec:intro}
Collisions are an essential feature of celestial mechanics problems,
and the scourge of global dynamics. This is because, 
given the initial positions and velocities of a collection of $N$ gravitating 
point masses, collisions obstruct the extension of 
local solution curves throughout all time.  It is then
a fundamental question to determine the embedding of the
forward collision manifold: 
that is, the set of all initial conditions which reach collision
in finite forward time.  The backward collision manifold, or ejection 
manifold, is defined similarly by considering backwards time evolution.

Despite intense of work on the problem over nothing less than centuries, 
the collision manifold is completely understood in only the simplest cases. For example 
the case of two bodies, which reduces to the Kepler problem, and a few other completely 
integrable systems.  An illuminating discussion of the state-of-the art for $N = 3$ bodies
is found in the introduction of the paper by Guardia, Kaloshin, and Zhang \cite{MR3951693},
where they prove a kind of density result for the collision manifold in a 
perturbative regime where two of the three bodies are very small.

One fruitful approach which provides valuable partial results 
is to focus on the complement of the collision manifold, 
and to study the properties of orbits whose global existence is
assumed a-priori. A capstone example of this strategy
is the celebrated classification theorem of 
Chazy, published in 1922 \cite{MR1509241}, describing the possible forward/backward 
asymptotic time behavior of all collision-free orbits in the three body problem.  
See again \cite{MR3951693} for a precise statement of Chazy's 
theorem in modern language, and for a thorough discussion of the surrounding literature.

Another, more constructive approach is to to study invariant sets like relative 
equilibria, relative periodic orbits, relative invariant tori, 
and heteroclinic/homoclinic connections between these 
(\textit{relative} invariance here refers to solutions studied in an appropriate rotating frame). 
Moreover, given one of these invariant objects it is typically possible to develop
a normal form describing the dynamics nearby.  
This approach has deep roots in Poincar\'{e}'s groundbreaking 
 New Methods \cite{MR1194622,MR1194623,MR1194624}.
From the point of view of collisions, orbits lying in fully invariant sets are
necessarily collision-free, and hence in the complement of both the ejection
and collision manifolds.

An antipodal, but equally constructive approach is to study the (local) collision manifold 
using techniques of regularization. The idea of regularizing 
three body collisions goes back (at least)
to the work of Levi-Civita in the 1920s \cite{MR1555161}, 
and we refer to Chapter three Szebehely's book \cite{theoryOfOrbits},
and to the lecture notes of Celletti \cite{10.1007/978-1-4020-4706-0_7}
for more on the historical development of regularization and many additional references.
The main point, from the perspective of the present discussion, is that 
regularizing coordinate transformations replace collision sets 
with well defined geometric objects which can be advected or ``grown'' using 
the regularized flow.  Further away from these regularized geometric objects, 
the local representation can be transformed back to the original coordinates and advected 
further.  In this light, collision manifolds are similar to 
local stable/unstable manifolds attached to invariant objects,
and exhibit all the complexity  
which Poincar\'{e} himself complained was ``difficult to draw''.

Rather than study the full complexity of the ejection and
collision manifolds, one can mimic the approach of dynamics systems 
theory and study only their intersections.
These intersections give rise to ``connecting orbits'' 
which we refer to as ejection-collision orbits.
Ejection-collision orbits where the same 
two bodies collide in forward and backward time are referred to as ``homoclinic'',  
and have been studies extensively from both perturbative and 
numerical view points.  We refer the interested reader to  the 
works of Oll\'e, Rodr\'{i}guez, and Soler
\cite{MR4162341,MR4110029} and to Seara, Oll\'{e}, Rodr\'{i}guez, and Soler
\cite{MR4518121} for powerful existence results, numerical 
simulations, and detailed discussion of the literature.  
See also the introduction of \cite{MR4576879}, 
discussed further below.

Another, more global situation is when a body $A$ collides with a body $B$ in backward time 
but with a different body $C$ in forward time.  We refer to such orbits as 
ejection-collision heteroclinics, and they are the focus of the present work.  

Analytical results for ejection-collision heteroclinics in the three body problem 
are found in Llibre, and Llibre and Lacomba \cite{MR0682839,MR0949626}.
Their arguments are perturbative,
with the rotating Kepler problem serving as the unperturbed system,
and the masses of the second and third bodies
along with the reciprocal of energy as small parameters. 
In lay terms, their result shows that an energetic enough particle 
aimed at the Moon from the Earth can strike it.

For larger mass ratios and/or lower energies, purely analytic techniques falter.  
Numerical computations come to the rescue and, 
when mathematical rigor is desired, 
techniques from validated numerics (for example interval arithmetic)
can be combined with a-posteriori, Newton-Kantorovich type analysis to construct 
existence arguments.  This is a part of what is called 
\textit{computer-assisted proof in analysis}, and some remarks about 
the literature are given in Appendix \ref{sec:literature}.
In the present work, we develop computer assisted existence 
results for parameter dependent families of 
heteroclinic ejection-collision orbits in the planar circular 
restricted three body problem (CRTBP).

\begin{figure}[t]
\begin{center}
\includegraphics[height=5.0cm]{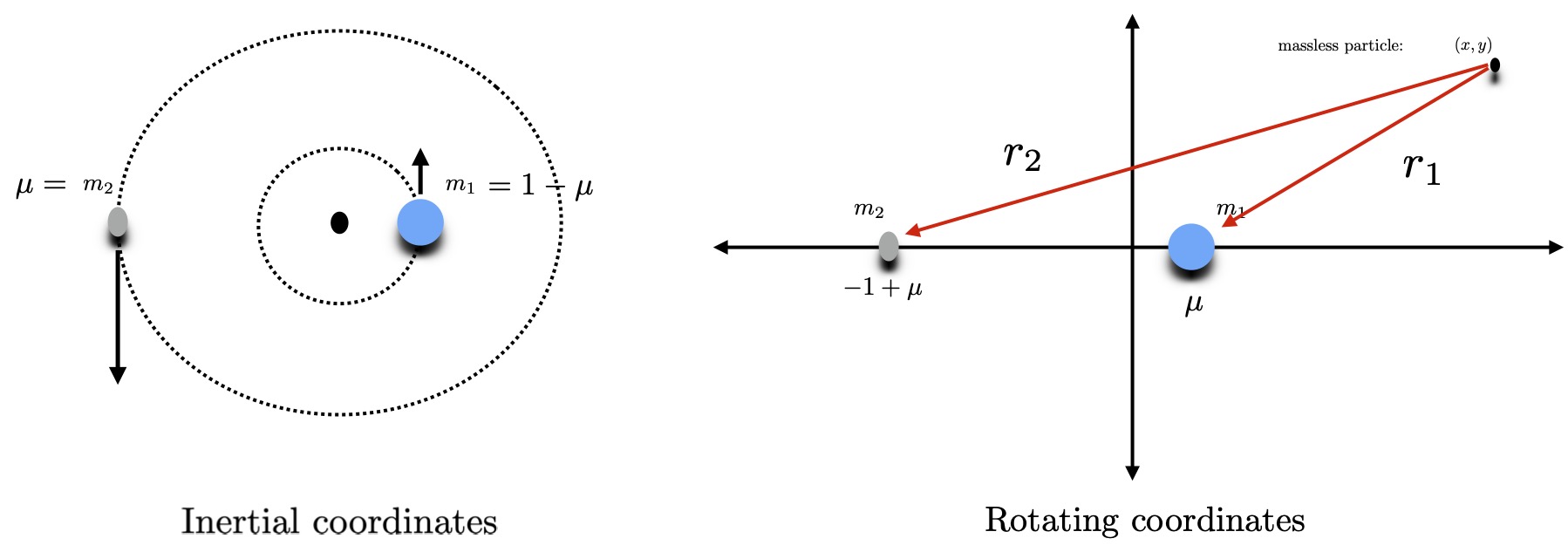}
\end{center}
\caption{\textbf{Configuration of the Circular Restricted Three Body Problem: }
the left frame illustrates the motion of the massive primary bodies $m_1$ and 
$m_2$ in inertial coordinates.  The are restricted a-priori to Keplerean circles.
The line determined by the position of 
$m_1$, $m_2$ and the center of mass rotates at a constant speed, 
and it is possible to change to a co-rotating coordinate system with 
constant angular velocity, thus removing the motion of primaries. 
Inserting a massless test particle into the resulting gravitational 
vortex, one obtains the CRTBP as illustrated in the right frame. }
\label{fig:CRTBP}
\end{figure}

\begin{figure}[t]
\begin{center}
\includegraphics[height=4.0cm]{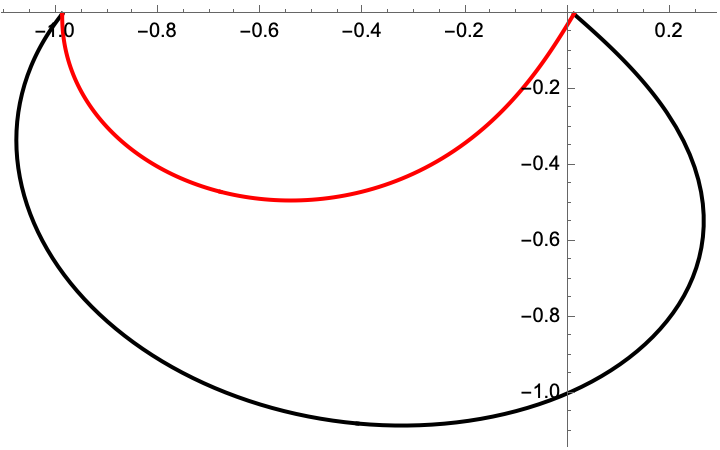}
\includegraphics[height=4.0cm]{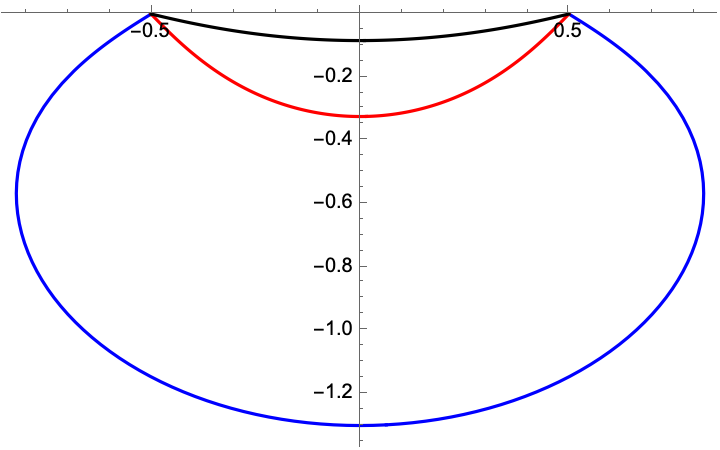}
\end{center}
\caption{\textbf{Branches parameterized by energy: }
The 5 approximate heteroclinic ejection-collision
orbits form the statement of Theorem \ref{th:other}. 
The left frame illustrates trajectories 1 (black) and 2 (red) in 
the Earth-Moon system, while the 
right frame illustrates trajectories 3 (black), 4 (red), and 5 (blue)
in the equal mass (or Copenhagen) system.}
\label{fig:sol345}
\end{figure}

\begin{figure}[t]
\begin{center}
\includegraphics[height=5.0cm]{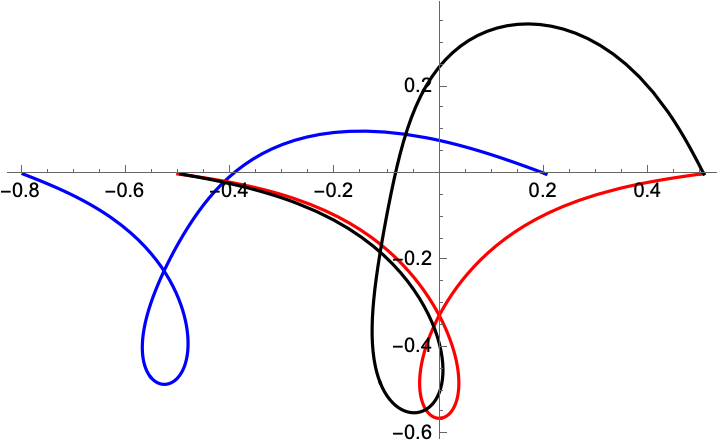}
\end{center}
\caption{\textbf{A global branch of heteroclinic ejection-collision orbits parameterized by 
mass ratio:}
The figure illustrates some approximate orbits from Theorem \ref{th:fixedC}. 
We note that, in the CRTBP, the mass parameter $\mu$ is normalized so that it
lies in the interval $[0, 1/2]$.  This branch starts a $\mu = 1/2$, decreases to $\mu = \mu_* < 1/2$, 
undergoes a fold bifurcation, and returns to $\mu = 1/2$: hence we follow it over its full 
lifespan.  The blue curve illustrates 
the solution at the bifurcation parameter (mass ratio) $\mu=\mu^*$, while
the red and black the solutions at $\mu=1/2$.
The figure illustrates the result of Theorem \ref{th:fixedC}, which shows that there 
exists a one parameter family 
of ejection-collision orbits, parameterized by $\mu$, which  
connects the red to the black solution through the saddle node bifurcation at the blue solution.  Note that 
for the orbits in this family, the radial velocity with respect to the large primary vanishes twice, hence the 
family is dynamically different from the one established in  \cite{MR0682839}.} 
\label{fig:bb}
\end{figure}

The CRTBP was introduced by Poincar\'{e} as a simplified model of three body 
dynamics, and we recall the equations of motion in Section \ref{sec:PCRTBP}.
Briefly, the problem studies --in a co-rotating frame --  the dynamics 
of an infinitesimal particle moving in the vicinity of two massive 
primary bodies (think of a man-made satellite or comet) .  
The primary bodies are assumed to evolve on Keplerian circles, and 
to be unaffected by the presence of the infinitesimal third body.
Again, an excellent general reference remains the book of 
Szebehely \cite{theoryOfOrbits}, with a more modern treatment 
being the book of Meryer, Hall, and Offin \cite{MR2468466}.
See Figure \ref{fig:CRTBP} for a schematic illustration of the 
CRTBP.  
Example results obtained using our methods 
are illustrated in Figures  \ref{fig:sol345} and \ref{fig:bb}, with 
detailed descriptions found in Section \ref{sec:results}.

We note that the the recent work of Capi\'{n}ski, Kepley, and the second
deals with computer assisted proofs for collision and near collision orbits in 
the CRTBP with energy and mass parameters fixed \cite{MR4576879}.  
Additional comments about the differences between \cite{MR4576879}
and the present work are found in Appendix \ref{sec:literature}. 
We also mention a recent preprint by Capi\'{n}ski and Pasiut \cite{MaciejAlexander} 
which applies techniques based on those of \cite{MR4576879} to prove the 
existence of a horse shoe-like invariant set where orbits
 pass arbitrarily close to collision infinitely many times
 in the Earth-Moon CRTBP .

Finally we note that, because our computer assisted 
arguments are built on top of accurate numerical approximations 
of one parameter branches of ejection-collision orbits in the CRTBP,
our work is closely related (and deeply indebted)
to classical numerical continuation and bifurcation theory:
especially to the literature on computation and continuation of invariant 
objects in celestial mechanics.  This literature in this area is 
vast, going all the way back to the work of Darwin  \cite{MR1554890}, 
Str\"{o}mgren \cite{stromgrenRef}, and Moulton \cite{moultonBook}
at the turn of the Twentieth Century.  
A thorough review of this literature 
is a task far beyond the scope of the present work, and we 
refer to the works of 
Mu\~{n}oz-Almaraz, Freire, Gal\'{a}n, Doedel, and Vanderbauwhede,
\cite{MR2003792}, 
Doedel, Romanov, Paffenroth, Keller, Dichmann, Gal\'an-Vioque, and Vanderbauwhede
\cite{MR2355523}, 
Doedel, Paffenroth, Keller, Dichmann, Gal\'an-Vioque, and Vanderbauwhede
\cite{MR1992054}, and 
Calleja, Doedel, Humphries, Lemus-Rodr\'{\i}guez, A. and Oldeman \cite{MR2989589}
for detailed descriptions of numerical continuation techniques for 
branches of invariant objects in the CRTBP.  By consulting these 
 papers, and the references therein, the reader will obtain an much 
clearer picture of the wider literature.    
More general references include the Lecture notes of Simo
\cite{1990mmcm.conf..285S} and of Doedel \cite{doedel_paris}, as well 
as the book of Kuznetsov \cite{MR4592585}.

\bigskip

\bigskip

\begin{figure}[t!]
\begin{center}
\includegraphics[height=7.5cm]{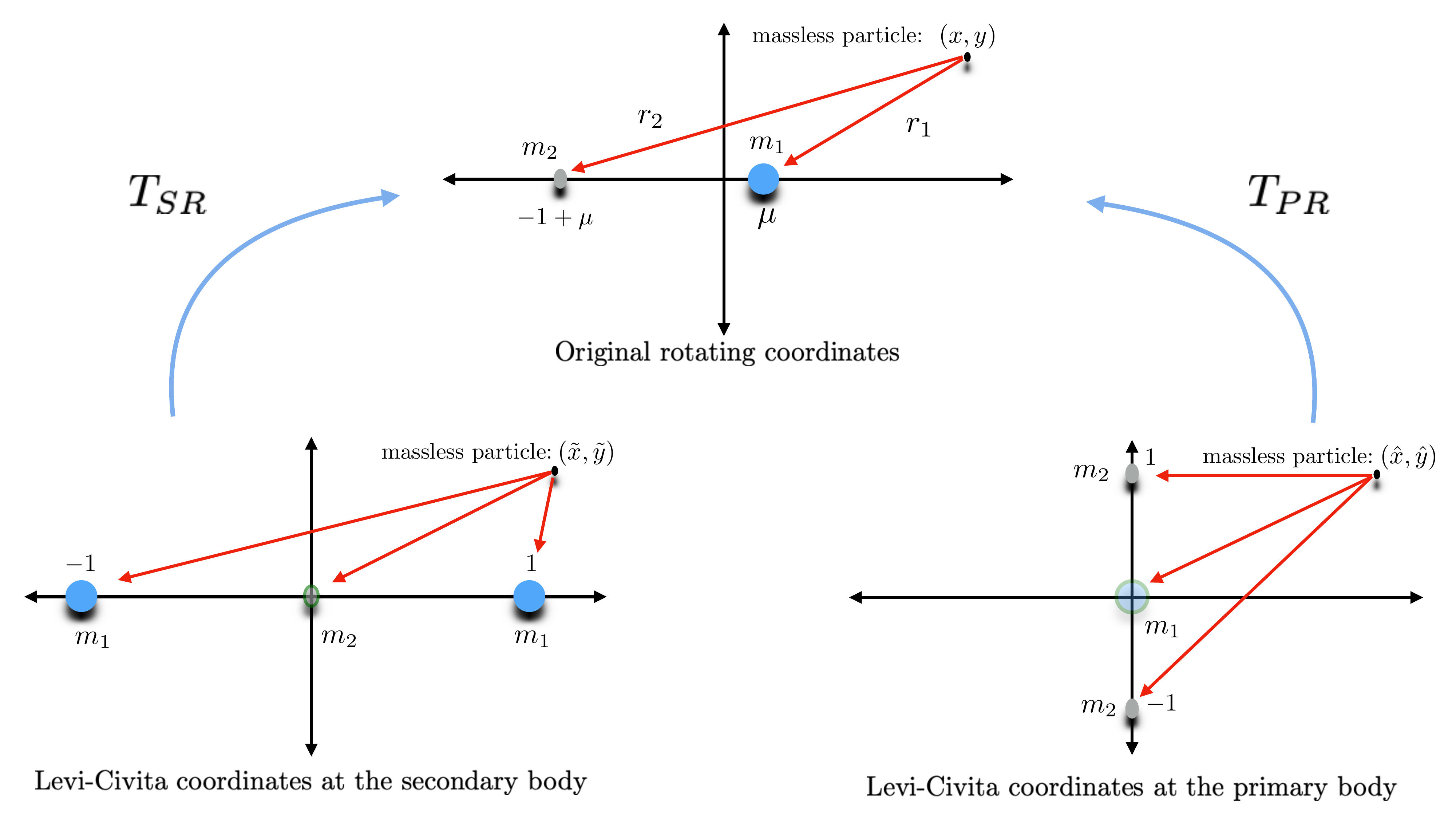}
\end{center}
\caption{\textbf{Levi-Civita Regularized Coordinates}: center top
figure depicts the configuration space of the CRTBP, with singularities
at $(x,y) = (-1 + \mu, 0)$ and at $(x,y) = (\mu, 0)$ due to collision of the 
massless particle with the smaller primary $m_2$ or the 
larger primary $m_1$ respectively.  The bottom left  and right figures 
depict the configuration spaces of the Levi-Civita systems associated
with regularization of the smaller and larger masses respectively.
Note that in the bottom left and  right frames the massive bodies 
$m_2$ and $m_1$ respectively have been moved to the origin.
These bodies are depicted as ``ghosts'' using transparencies to 
indicate that the resulting fields are perfectly well defined there.
However, the double covering introduced by the complex squaring 
function creates mirror images of the remaining primary, so that 
the Levi-Civita systems still have a pair of singular point.  
These mirrored singular points are located at $(x,y) = (\pm 1, 0)$ for the 
case the regularization at $m_2$ and at $(x,y) = (0, \pm 1)$ for the 
regularization at $m_1$.  The mirrored singularities play no role in our set up, 
as we employ the Levi-Civita coordinates only in a small neighborhood of the 
origin.}
\label{fig:PCRTBP_coordinates}
\end{figure}

\smallskip              
                
The remainder of the paper is organized as follows.
In Section \ref{sec:PCRTBP} we recall the equations of 
motion for the circular restricted three body problem,
as well as the coordinate changes and the resulting 
regularized equations of motion resulting from the
Levi-Civita transformations.  
Section \ref{sec:results} states precisely our main results, 
and in Section \ref{sec:BVP} we describe the formulation of 
the fixed point problem for the boundary value problems
(BVPs) used to study ejection-collision orbits.
In Section \ref{sec:bif} we present the 
appropriate bifurcation theory for branches of solutions to the 
BVPs discussed in Section \ref{sec:BVP}, and in Section 
\ref{sec:branches} we discuss the a-posteriori analysis for branches
of such solutions.  Finally, in Section \ref{sec:cap} we discuss some 
implementation details for the computer assisted proofs.  In particular, 
we describe computer assisted verification of the hypotheses of the 
a-posteriori theorems developed earlier in the paper.   
Appendix \ref{sec:literature} provides a few additional remarks 
about the surrounding literature.

\section{The planar circular restricted three body problem} \label{sec:PCRTBP}

Let $m_1$ and $m_2$ denote the masses of the primary bodies.
We assume that $m_1 \geq m_2$, and define the mass ratio
\[
\mu = \frac{m_2}{m_1 + m_2}.
\]
Note that $\mu \in [0, 1/2]$, with $\mu = 0$ when $m_2 = 0$
and $\mu = 1/2$ when $m_1 = m_2$. 

Consider
\begin{equation} \label{eq:CRTBP_ODE}
\gamma' = f_\mu(\gamma),
\end{equation}
where $\gamma = (x,p,y,q)$, and $f_\mu$ is the 
one parameter family of vector fields 
given by 
\begin{equation}  \label{eq:PCRTBP}
f_\mu(x,p,y,q):=\left( 
\begin{array}{c}
p \\ 
2q+x-\frac{(1-\mu )\left( x-\mu \right) }{r_1^3}-\frac{%
\mu \left( x+1-\mu \right) }{r_2^3} \\ 
q \\ 
-2p+y-\frac{(1-\mu )y}{r_1^3}-\frac{\mu y}{r_2^3}%
\end{array} %
\right).
\end{equation}
Here
\begin{equation*}
r_{1}=\sqrt{(x-\mu )^{2}+y^{2}},\quad \quad \text{and}\quad \quad
r_{2}=\sqrt{(x+1-\mu )^{2}+y^{2}}.
\end{equation*}
We refer to the ordinary differential equation 
given by Equation \eqref{eq:CRTBP_ODE}
as the planar circular restricted three body 
problem or CRTBP, and remark that when $\mu = 0$ the 
system reduces to the rotating Kepler problem.
This is the key to perturbative approaches to the problem, 
but will play no role in the present work.

Equation \eqref{eq:CRTBP_ODE} is derived, starting from Newton's 
laws for three massive gravitating bodies, 
by assuming the motion of the primary bodies is known and circular, and that 
the mass of the third body goes to zero.
One then changes to a rotating frame revolving with the frequency 
of the primaries based at their center of mass.  
Equation \eqref{eq:CRTBP_ODE} is obtained after
choosing units of distance, mass, and time so that 
the primaries (whose locations are fixed in by the rotating frame)
have coordinates $(x_1, y_1) = (\mu, 0)$ and $(x_2, y_2) = (-1+\mu, 0)$,
and so that the masses become 
$m_1 = 1- \mu$ and $m_2= \mu$.
A schematic illustration of the configuration space is 
seen in the right frame of Figure \ref{fig:CRTBP}.

Orbits of the massless particle -- that is, solutions of the differential equation
\eqref{eq:CRTBP_ODE} --
conserve the scalar quantity 
\begin{equation}
E\left(x,p,y,q\right) =-p^{2}-q^{2}+2\Omega (x,y),
\label{eq:JacobiIntegral}
\end{equation}
where
\begin{equation*}
\Omega (x,y)=(1-\mu )\left( \frac{r_{1}^{2}}{2}+\frac{1}{r_{1}}\right) +\mu
\left( \frac{r_{2}^{2}}{2}+\frac{1}{r_{2}}\right).
\end{equation*}%
The function $E$ defined in Equation \eqref{eq:JacobiIntegral}
is referred to as the Jacobi integral or Jacobi constant and 
sometimes, in a slight abuse of terminology, we refer to $E$
as the energy of the system.

We now define the two dimensional affine subsets of $\mathbb{R}^4$ given by 
\[
\mathcal{C}_\mu^1 = \left\{ (x,p,y,q) \in \mathbb{R}^4 \, | \, x = \mu \mbox{ and } y = 0 \right\}, 
\]
and
\[
\mathcal{C}_\mu^2 = \left\{ (x,p,y,q) \in \mathbb{R}^4 \, | \, x =-1+\mu \mbox{ and } y = 0 \right\}.
\]
These are the singular sets associated with the positions of the
first and second primary bodies respectively. They correspond to 
the situation where the position of the 
massless particle coincides with the position of one of the primaries:
that is, a collision.  

Let 
\[
\mathcal{C}_\mu = \mathcal{C}_\mu^1 \cup \mathcal{C}_\mu^2,
\]
and note that 
\begin{equation} \label{eq:domain}
U_\mu = \mathbb{R}^4 \backslash \mathcal{C}_\mu,
\end{equation}
is the natural domain of $f_\mu$. In this paper we focus on 
$\mu \in (0, 1/2]$, as the case of 
$m_2 = \mu = 0$ degenerates to the rotating Kepler problem, 
and the singularity associated with the second primary vanishes.
In this case no ejection from or collision with the 
second primary is possible.

We note that, since $f_\mu$ is locally Lipschitz on $U_\mu$
 (in fact real analytic), we have that the orbit of every 
initial condition $\gamma(0) \in U_\mu$
either exists for all time,
or else accumulates to $\mathcal{C}_\mu$ in finite time.  
 We refer to  $\mathcal{C}_\mu$ as \textit{the collision set}, and 
 say that an orbit which accumulates to $\mathcal{C}_\mu$
 in finite forward time is a \textit{collision orbit}.  
Similarly, an orbit which accumulates to $\mathcal{C}_\mu$
 in finite backward time is referred to as an \textit{ejection
 orbit}.  Finally, an orbit which accumulates to $\mathcal{C}_\mu$
 in both finite forward and finite backward time is 
 an \textit{ejection-collision orbit}, and the  
 objects of our study.

To formalize these notions, let $\pi \colon \mathbb{R}^4 \to \mathbb{R}^2$
be the projection onto position variables given by 
\[
\pi (x, p, y, q) = (x,y).
\]
We make the following definition.

\begin{definition}[Ejection-collisions]\label{def:EC}
Let $T_1 \in (-\infty, 0)$, $T_2 \in (0, \infty)$
 and $\mu\in(0,1/2]$.
We say that the curve 
$\gamma\colon (T_1, T_2)\to U_\mu \subset \mathbb{R}^4$
is an $m_1$ to $m_2$ ejection-collision orbit if 
\[
\frac{d}{dt} \gamma(t) = f_\mu(\gamma(t))\,,
\]
for $t \in (T_1, T_2)$, and 
\[
\lim_{t \to T_1} \pi \gamma(t) = (\mu, 0), 
\quad \quad \mbox{and} \quad \quad
\lim_{t \to T_2} \pi \gamma(t) = (-1+\mu, 0).
\]
The energy of the ejection-collision orbit is 
\[
E(\gamma(0)) = C.
\]
If the limits are reversed then $\gamma$ is an $m_2$ to $m_1$ ejection-collision orbit.  
\end{definition}

\medskip

To understand velocity as 
the infinitesimal body approaches collision, 
consider $\gamma \colon (T_1, T_2) \to U_\mu \subset \mathbb{R}^4$ an
$m_1$ to $m_2$ ejection-collision orbit in the sense of Definition \ref{def:EC}.
We claim that a collision is in fact just a certain kind of finite time blow up.  
To see this, let $v(t) = \sqrt{p(t)^2 + q(t)^2}$ denote the magnitude of velocity.
(The discussion is easily modified for $m_2$ to $m_1$ ejection-collisions).
The energy of the orbit is   
$C = E(\gamma(0))$, and we write 
\[
\gamma(t) = \left(
\begin{array}{c}
x(t) \\
p(t) \\
y(t) \\
q(t)
\end{array}
\right),
\]
to denote the components of the orbit.
In particular, since $\gamma$ is an $m_1$ to $m_2$ ejection-collision we have that 
\[
\lim_{t \to T_1} x(t) = -1 + \mu, 
\quad \quad \mbox{and} \quad \quad 
\lim_{t \to T_2} x(t) = -\mu, 
\]
so that 
\begin{equation} \label{eq:velLimits}
\lim_{t \to T_1} r_1(t) = 
\lim_{t \to T_2} r_2(t) = 0.
\end{equation}
(These are swapped if $\gamma$ is $m_2$ to $m_1$).

Recalling Equation \eqref{eq:JacobiIntegral}, we now have that 
\[
C =-v(t)^2+2\Omega (x(t),y(t)),
\]
where
\[
\Omega (x(t),y(t))=(1-\mu )\left( \frac{r_{1}(t)^{2}}{2}+\frac{1}{r_{1}(t)}\right) +\mu
\left( \frac{r_{2}(t)^{2}}{2}+\frac{1}{r_{2}(t)}\right).
\]
It follows from Equation \eqref{eq:velLimits} that 
 $\Omega(x(t), y(t)) \to \infty$ as $t \to T_{1,2}$, and since 
\[
v(t)^2 =2\Omega (x(t),y(t)) - C, 
\] 
with $C$ constant, we have that 
\[
v(t) \to \pm \infty, \quad \quad \quad \mbox{as } t \to T_{1,2}.
\]
That is, at least one of $p(t)$ or $q(t)$ becomes infinite as 
the infinitesimal body approaches collision.  One concludes that
orbits accumulate to $C_\mu^{1,2}$ only at infinity.

\subsection{Regularized coordinates for collisions with $m_{1}$}
\label{sec:regularized1}
We now review the classical Levi-Civita transformations
\cite{MR1555161}, which allow us to extend orbits of the CRTBP
up to and through collision.  We state a number of standard results 
without proof, and refer again to \cite{theoryOfOrbits} and
\cite{10.1007/978-1-4020-4706-0_7} for more details.

To regularize a collision with $m_{1}$, 
write $z=x+iy$ and define new regularized
variables $\hat{z}=\hat{x}+i\hat{y}$.
These are related to $z$ by the transformation 
\begin{equation*}
\hat{z}^{2}=z-\mu .
\end{equation*}%
It is also necessary to rescale time in the
new coordinates according to the formula
\begin{equation*}
\frac{dt}{d\hat{t}}=4|\hat{z}|^{2}.
\end{equation*}
Next one fixes a value $c$ of the Jacobi constant and, 
a standard calculation (see \cite{theoryOfOrbits}),
transforms the vector field $f$
of Equation \eqref{eq:PCRTBP}  to the regularized Levi-Civita vector
field $f_{1}^{c}\colon U_{1}\rightarrow \mathbb{R}^{4}$ given by
\[
f_1^c(\hat x, \hat p, \hat y, \hat q) = 
\]
\begin{eqnarray}  \label{eq:regularizedSystem_m1}
\left(
\begin{array}{c}
\hat{p},  \\
8\left( \hat{x}^{2}+\hat{y}^{2}\right) \hat{q}+12\hat{x}%
(\hat{x}^{2}+\hat{y}^{2})^{2}+16\mu \hat{x}^{3}+4(\mu -c)\hat{x}
+\frac{8\mu (\hat{x}^{3}-3\hat{x}\hat{y}^{2}+\hat{x})}{((\hat{x}^{2}
+\hat{y}^{2})^{2}+1+2(\hat{x}^{2}-\hat{y}^{2}))^{3/2}} \\
\hat{q} \\
-8\left( \hat{x}^{2}+\hat{y}^{2}\right) \hat{p}
+12\hat{y}\left( \hat{x}^{2}+\hat{y}^{2}\right) ^{2}-16\mu \hat{y}^{3}+4\left( \mu
-c\right) \hat{y} 
+\frac{8\mu (-\hat{y}^{3}+3\hat{x}^{2}\hat{y}+\hat{y})}{((\hat{x}^{2}+\hat{%
y}^{2})^{2}+1+2(\hat{x}^{2}-\hat{y}^{2}))^{3/2}}
\end{array}
\right). 
\end{eqnarray}
Note that the open set $U_{1}\in \mathbb{R}^{4}$ defined by 
\begin{equation*}
U_{1}=\left\{ \mathbf{\hat{x}}=(\hat{x},\hat{p},\hat{y},\hat{q})\in \mathbb{R%
}^{4} : \left( \hat{x},\hat{y}\right) \notin \left\{ \left( 0,-1\right)
,\left( 0,1\right) \right\} \right\},
\end{equation*}%
serves as the domain of the regularized system, 
 that the regularized vector field is well defined at
the origin $\left( \hat{x},\hat{y}\right) =\left( 0,0\right) $, and that the origin is 
mapped to the collision with $m_{1}$ when inverting the Levi-Civita
coordinate transformation.
A schematic illustration of the regularized configuration space is given in the 
bottom right frame of Figure \ref{fig:PCRTBP_coordinates}.

By working out  the change of coordinates for the velocity variables, 
one defines the coordinate 
transformation $T_{1}\colon U_{1}\backslash \mathcal{C}_{1}\rightarrow U$
between the original and regularized coordinates by 
\begin{equation}
\mathbf{x}=T_{1}(\mathbf{\hat{x}}):=\left( 
\begin{array}{c}
\hat{x}^{2}-\hat{y}^{2}+\mu  \\ 
\frac{\hat{x}\hat{p}-\hat{y}\hat{q}}{2(\hat{x}^{2}+\hat{y}^{2})} \\ 
2\hat{x}\hat{y} \\ 
\frac{\hat{y}\hat{p}+\hat{x}\hat{q}}{2(\hat{x}^{2}+\hat{y}^{2})}%
\end{array}%
\right).  \label{eq:T1-def}
\end{equation}%
Note that $T_1$ is a local diffeomorphism on $U_{1}\backslash \mathcal{C}_{1}$. 
The regularized vector field conserves the first integral $E_{1}^{c}\colon
U_{1}\rightarrow \mathbb{R}$ given by 
\begin{eqnarray}
E_{1}^{c}(\mathbf{\hat{x}}) &=&-\hat{q}^{2}-\hat{p}^{2}+4(\hat{x}^{2}+\hat{y}%
^{2})^{3}+8\mu (\hat{x}^{4}-\hat{y}^{4})+4(\mu -c)(\hat{x}^{2}+\hat{y}^{2}) 
\notag \\
&&+8(1-\mu )+8\mu \frac{(\hat{x}^{2}+\hat{y}^{2})}{\sqrt{(\hat{x}^{2}+\hat{y}%
^{2})^{2}+1+2(\hat{x}^{2}-\hat{y}^{2})}}.  \label{eq:reg_P_energy}
\end{eqnarray}

An essential point is that the fixed energy parameter $c$ 
appears in $f_{1}^{c}$ and $E_{1}^{c}$ as a parameter.  We remark that 
a point with energy $c$ in the original coordinate system is transformed to 
a point with energy zero in the regularized coordinates. 
Then, since collision in the original coordinates corresponds to the origin in the 
regularized coordinates, the intersection of the collision set with the zero 
energy level set is given by 
\[
E_{1}^{c}(0, \hat p, 0, \hat q) =  
-\hat{q}^{2}-\hat{p}^{2}
+8(1-\mu ) = 0,
\]
or 
\begin{equation} \label{eq:circ1}
\hat{q}^{2} + \hat{p}^{2} = 
8(1-\mu ).
\end{equation}
That is the the Levi-Civita 
transformation maps the intersection of the collision set 
 $\mathcal{C}_\mu^1$ with the energy level set $E = c$ to the 
 origin cross the velocity circle given in Equation \eqref{eq:circ1}.

\subsection{Regularized coordinates for collisions with $m_{2}$} \label{sec:regSecondPrimary}
\label{sec:regularized2}

In a precisely analogous way, 
one can regularize collisions with the second primary as follows.  Write 
$\tilde{z}=\tilde{x}+i\tilde{y}$, define $\tilde{z}^{2}=z+1-\mu $, rescale time
as $dt/d\tilde{t}=4|\tilde{z}|^{2}$, and
fix an energy level $c$ in the original coordinates.  Another 
lengthy calculation results in the 
regularized field $f_{2}^{c}:U_{2}\rightarrow \mathbb{%
R}^{4}$  defined as
\[
f_2^c(\tilde x, \tilde p, \tilde y, \tilde q) = 
\]
\begin{eqnarray}  \label{eq:regularizedSystem_m2}
\left(
\begin{array}{c}
\tilde{p}
\\
8\left( \tilde{x}^{2}+\tilde{y}^{2}\right) \tilde{q}
+12\tilde{x}(\tilde{x}^{2}+\tilde{y}^{2})^{2}-16(1-\mu )\tilde{x}%
^{3}+4\left( (1-\mu )-c\right) \tilde{x} 
+\frac{8(1-\mu )\left( -\tilde{x}^{3}+3\tilde{x}\tilde{y}^{2}+\tilde{x}%
\right) }{((\tilde{x}^{2}+\tilde{y}^{2})^{2}+1+2(\tilde{y}^{2}-\tilde{x}%
^{2}))^{3/2}}  \\
\tilde{q}   \\
-8\left( \tilde{u}^{2}+\tilde{y}^{2}\right) \tilde{p}%
+12\tilde{y}(\tilde{x}^{2}+\tilde{y}^{2})^{2}+16(1-\mu )\tilde{y}%
^{3}+4\left( (1-\mu )-c\right) \tilde{y}  
+\frac{8(1-\mu )\left( \tilde{y}^{3}-3\tilde{x}^{2}\tilde{y}+\tilde{y}%
\right) }{((\tilde{x}^{2}+\tilde{y}^{2})^{2}+1+2(\tilde{y}^{2}-\tilde{x}%
^{2}))^{3/2}}
\end{array}
\right)
\end{eqnarray}%
defined on the domain 
\begin{eqnarray*}
U_{2}&:= & \mathbb{R}^4 \backslash \mathcal{C}_2 = \left\{ \mathbf{\tilde{x}}=(\tilde{x},\tilde{p},\tilde{y},\tilde{q}%
)\in \mathbb{R}^{4}\,|\,\left( \tilde{x},\tilde{y}\right) \notin \left\{
\left( -1,0\right) ,\left( 1,0\right) \right\} \right\} , \quad \quad \mbox{where} \\
\mathcal{C}_{2}&:= &\left\{ \mathbf{\tilde{x}}=(\tilde{x},\tilde{p},\tilde{y}%
,\tilde{q})\in \mathbb{R}^{4}\,|\,\tilde{x}=\tilde{y}=0\right\}.
\end{eqnarray*}%
The complete
change of coordinates is given by 
\begin{equation}
\mathbf{x}=T_{2}\left( \mathbf{\tilde{x}}\right) =\left( 
\begin{array}{c}
\tilde{x}^{2}-\tilde{y}^{2}+\mu-1 \\ 
\frac{\tilde{x}\tilde{p}-\tilde{y}\tilde{q}}{2(\tilde{x}^{2}+\tilde{y}^{2})}
\\ 
2\tilde{x}\tilde{y} \\ 
\frac{\tilde{y}\tilde{p}+\tilde{x}\tilde{q}}{2(\tilde{x}^{2}+\tilde{y}^{2})}%
\end{array}%
\right)  \label{eq:T2-def}
\end{equation}
A schematic illustration of the configuration space is given in the 
bottom left frame of Figure \ref{fig:PCRTBP_coordinates}.

We also remark that system conserves the scalar quantity 
\begin{align}
E_{2}^{c}\left( \mathbf{\tilde{x}}\right) & =-\tilde{p}^{2}-\tilde{q}^{2}+4(%
\tilde{x}^{2}+\tilde{y}^{2})^{3}+8(1-\mu )(\tilde{y}^{4}-\tilde{x}%
^{4})+4\left( (1-\mu )-c\right) (\tilde{x}^{2}+\tilde{y}^{2})  \notag \\
& \quad +8(1-\mu )\frac{\tilde{x}^{2}+\tilde{y}^{2}}{\sqrt{(\tilde{x}^{2}+%
\tilde{y}^{2})^{2}+1+2(\tilde{y}^{2}-\tilde{x}^{2})}}+8\mu,  \label{eq:E2}
\end{align}
so that, arguing as in the previous section, the intersection of the 
collision set $\mathcal{C}_\mu^2$ with the $c$ energy level set 
is given by the product of the origin and the circle 
\begin{equation}
\tilde p^2 + \tilde q^2 = 8 \mu.
\end{equation}

\section{Statement of the main results} \label{sec:results}

\medskip

\begin{table}[h!]
  \begin{center}
    \label{tab:table1}
    \begin{tabular}{c|c|c} 
          Solution & $\mu$ & $C$\\
      \hline
      1 & $\mu_{em}$ & 2\\
      2 & $\mu_{em}$ & $845\times2^{-9}\approx1.65$\\
      3 & 1/2 & -6\\
      4 & 1/2 & 2\\
      5 & 1/2 & 2\\
    \end{tabular}
  \end{center}
\end{table}

We now describe our main results.  
The goal of the remainder of present work is to prove the following theorems.

\begin{theorem}\label{th:other}
For each pair of parameter values $(\mu,C)$ given in
Table \ref{tab:table1}, there exists an $m_1 = 1-\mu$ to $m_2 = \mu$ 
ejection-collision orbit at energy level $C$.
Solutions 3,4,5 are symmetric with respect to $x\mapsto -x$, see Figure \ref{fig:sol345}.
\end{theorem}

\medskip

Having established the existence of ejection-collision orbits as illustrated in Theorem \ref{th:other},
our next task is to ask wether such orbits persist as the energy or mass ratio is varied.
The next three theorems, which constitute the main results of the present work, 
deal with branches of ejection-collision orbits from one primary to the other, and bifurcations 
of such branches.  We consider two kinds of branches: those parameterized by energy with 
mass ratio held fixed, and those parameterized by mass ratio with energy fixed.  
For example we have the following.

\begin{theorem}\label{th:fixedC}
For $C_* = 205 \times 2^{-6} \approx 3.203$,
there exists 
\[
\mu_*  \in  [3261, 3263 ] \times2^{-14} \quad \quad (\mbox{that is }  \mu_* \approx  0.1991),
\]
and two distinct branches of $m_2$ to $m_1$ 
ejection-collision orbits parameterized by $\mu \in [\mu_*, 1/2]$
and energy $C_*$.
The two branches are real analytic on $(\mu_*, 1/2]$ and terminate in a saddle node bifurcation 
at $\mu_*$. 
\end{theorem}

The next theorem establishes the existence of
a one parameter branch  of ejection-collision orbits, parameterized by 
energy, in the Copenhagen problem.  That is, 
the CRTBP with $\mu = 1/2$ (equal masses).  

\begin{theorem}\label{th:muonehalf}
For $\mu = 1/2$, there exist 
\[
C_* = 2163799\times2^{-20}
\quad \quad (\mbox{that is } C_* \approx2.063559532), 
\]
and
\[
C^*  \in [8453, 8454] \times2^{-12}, \quad \quad
(\mbox{that is } C^*  \approx 2.06372),
\]
and two distinct 
branches of $m_1$ to $m_2$ ejection-collision orbits for $f_{1/2}$
parameterized by $C \in [C_*, C^*]$.
The branches are real analytic on $(C_*, C^*)$ and 
terminate in a saddle node bifurcation at $C^*$.
\end{theorem}

The last theorem is similar, but for the Earth/Moon problem. 
That is, the CRTBP with $\mu = \mu_{em}\defeq13256063\times2^{-30}\approx0.01234567$.  
                
\begin{theorem}\label{th:muearthmoon}
For $\mu_{em} \defeq13256063\times2^{-30}\approx0.01234567$, there exist
\[
C_* =34395443031\times2^{-34} \quad \quad (\mbox{that is } C_* \approx2.0020782849169336259), 
\]
and
\[
C^*  \in (34395443031\times2^{-34},34395443369\times2^{-34}),
\]
and two distinct 
branches of $m_1$ to $m_2$ ejection-collision orbits for $f_{ \mu_{em}}$
parameterized by $C \in [C_*, C^*]$.
The two branches are real analytic on $(C_*, C^*)$ and 
terminate in a saddle node bifurcation at $C^*$.
\end{theorem}

\begin{figure}[t]
\begin{center}
\includegraphics[height=5.0cm]{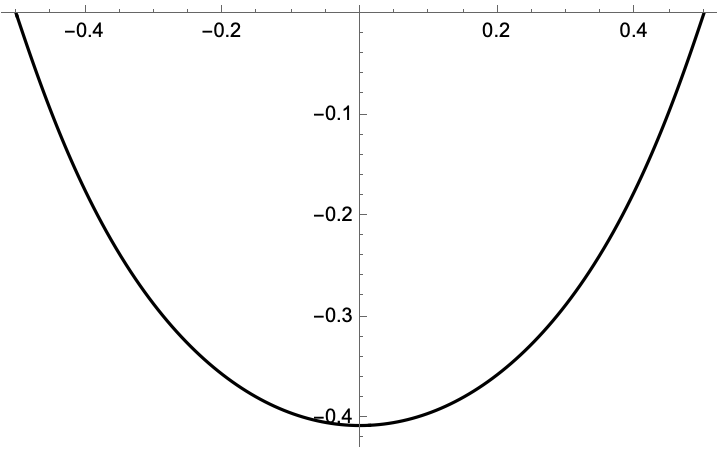}
\end{center}
\caption{Approximate orbits in the configuration space in rotating coordinates of the solution at $C^*$ obtained by Theorem \ref{th:muonehalf}.  For this family the radial velocity with respect to the 
first primary does not vanish, hence it is related to the family of   \cite{MR0682839}.  However, Theorem 
 \ref{th:muonehalf} holds at the highly non-perturbative energy of $\mu = 1/2$. }
\label{fig:sb1}
\end{figure}

\begin{figure}[t]
\begin{center}
\includegraphics[height=5.0cm]{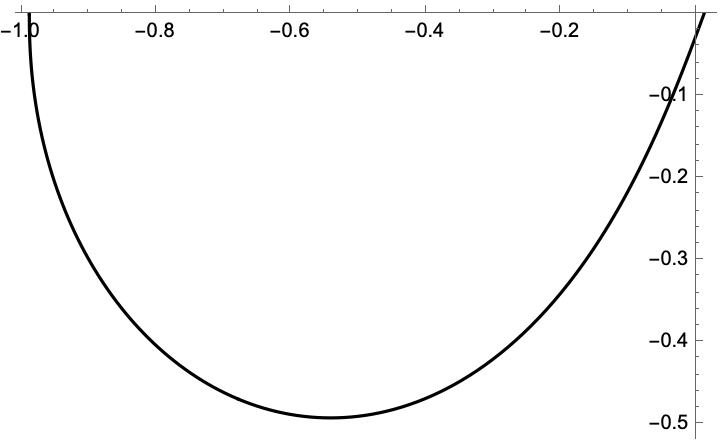}
\end{center}
\caption{Approximate orbits in the configuration space in rotating coordinates of the solution at $C^*$ obtained by Theorem \ref{th:muearthmoon}.  The family is similar to the one established by Theorem \ref{th:muonehalf}, 
but for the physically realistic Earth-Moon mass ratio.  }
\label{fig:sb2}
\end{figure}

\section{The fixed point equation} \label{sec:BVP}
In this Section we formulate a fixed point problem $\FF(X)=X$ whose solution
corresponds to an ejection-collision orbit in the CRTBP.
The main idea is to reformulate Definition \ref{def:EC}
in terms of the observations discussed in Section 
\ref{sec:PCRTBP}.  Again, the overall strategy is to approximately solve the fixed point 
equation using numerical methods, and then to construct a Newton-Kantorovich 
type argument in the vicinity of the approximate solution, with the goal of 
proving the existence of a true solution nearby.

Recall that the vector field defining CRTBP,
and both vector fields defining the regularized dynamics in Levi-Civita coordinates,
are real analytic on their domains of definition.  See Equations \eqref{eq:PCRTBP},  
\eqref{eq:regularizedSystem_m1}, and \eqref{eq:regularizedSystem_m2}.
Because of this, it is natural to look for real analytic solutions.  
So, for $\varrho,\rho>0$ define 
\begin{equation}\label{eq:tay}
\AA_\rho\defeq\{z:[-1,1]\to\mathbb{R}\,:\,z(t)=\sum_{j\ge0}z_j t^j\,,
\sum_{j\ge0}|z_j|\rho^j<+\infty\}\,,
\end{equation}
and
\begin{equation}\label{eq:cheb}
\BB_\varrho\defeq\{z:[-1,1]\to\mathbb{R}\,:\,z(t)=\sum_{j\ge0}z_j T_j(t)\,,
\sum_{j\ge0}|z_j|\varrho^j<+\infty\}.
\end{equation}
Here $\{z_j\}\subset\mathbb{R}$ and $T_j$ denotes the $j-$th \cheb polynomial.
These sets define Banach spaces (Banach algebras even) when 
endowed with the norms
$$
\|z\|_\AA\defeq\sum_{j\ge0}|z_j|\rho^j\,,\qquad\|z\|_\BB\defeq \sum_{j\ge0}|z_j|\varrho^j,
$$
and consist of real valued functions with analytic extensions to the complex unit disk of radius
$\rho$ in the case
of $\AA_\rho$, or to the complex ellipse with foci $\{-1,1\}$ and semiaxes $\half(\varrho+\varrho^{-1})$ 
where $\half(\varrho-\varrho^{-1})$ in the case of $\BB_\varrho$.

To study ejection-collision orbits from $m_1$ to $m_2$ we write 
$(x_1(t),p_1(t), y_1(t), q_1(t))$ to denote orbit segments in the 
Levi-Civita coordinates regularized
at $m_1$,  $(x_2(t),p_2(t), y_2(t), q_2(t))$ for orbit segments
in the coordinates regularized at $m_2$, and  
$(x_0(t),p_0(t), y_0(t), q_0(t))$, to denote orbit segments in the original coordinates.
We write $T_1$, $T_2$, and $T_0$ to denote the 
time spent in each these coordinate system (where $T_1$ and $T_2$ are
measured in regularized units of time).

We use Taylor series expansions to describe orbit segments in 
regularized coordinates, and Chebyshev series expansions in 
the original rotating coordinates.  This choice is motivated by the 
idea that we will take a small time step in the regularized coordinates
whose purpose is to move us away from collision, and then spend the remaining
time in the original rotating frame.  Then in general $T_0$ may be
substantially larger than $T_1, T_2$, and we use Chebyshev series 
to describe the longer orbit segment.
We rescale time so that $[0, T_1]$ and $[0, T_2]$ become
 $[0,1]$,  and $[0, T_0]$ becomes $[-1, 1]$, as these are natural 
 domains for Taylor and Chebyshev series expansions respectively.

The ejection-collision orbit is then 
a solution of the following constrained boundary value problem:
\begin{equation}\label{eq:first}
\begin{cases}
\dot p_1=T_1f_1(x_1,y_1,p_1,q_1)\\
\dot q_1=T_1g_1(x_1,y_1,p_1,q_1)\\
\dot p_0=T_0f_0(x_0,y_0,p_0,q_0)\\
\dot q_0=T_0g_0(x_0,y_0,p_0,q_0)\\
\dot p_2=T_2f_2(x_2,y_2,p_2,q_2)\\
\dot q_2=T_2g_2(x_2,y_2,p_2,q_2)\\
\dot x_i=T_ip_i\,,\quad i=0,1,2\\
\dot y_i=T_iq_i\,,\quad i=0,1,2\\
x_1(0)=0\\
 y_1(0)=0\\
  p_1(0)^2+q_1(0)^2=8(1-\mu)\,,\\
T_{PR}(x_1(1),p_1(1),y_1(1),q_1(1))=(x_0(-1),p_0(-1),y_0(-1),q_0(-1))\,,\\
T_{SR}(x_2(-1),p_2(-1),y_2(-1),q_2(-1))=(x_0(1),p_0(1),y_0(1),q_0(1))\,,\\
x_2(0)=0\\
 y_2(0)=0\\
 p_2(0)^2+q_2(0)^2=8\mu\,,
\end{cases}
\end{equation} 
where $f_{1,2},g_{1,2}$ are the second and fourth components of the functions $f^c_i$ defined 
by Equations \eqref{eq:regularizedSystem_m1} and \eqref{eq:regularizedSystem_m2}
respectively, and $f_0$, $g_0$ are the second and fourth components of the 
CRTB vector field defined in Equation \eqref{eq:PCRTBP}.
Note that the equations $p_1(0)^2+q_1(0)^2=8(1-\mu)$ and $p_2(0)^2+q_2(0)^2=8\mu$
 impose that the orbit segments
in regularized coordinates begin or end with collision at the correct Jacobi level, as discussed in 
 Sections \ref{sec:regularized1} and \ref{sec:regularized2}.  The transformations 
 $T_{PR}$ and $T_{SR}$ insure that the endpoints of the regularized 
 orbit segments match the end points of the segment in the standard rotating 
 coordinates.

Now, for $x\in\AA_\rho$ (or $\BB_\varrho$), we denote by
$D^{-1}_{c}x$ the inverse of the differentiation operator 
applied to $x$.  We normalize so that the zeroth order Taylor (resp. \chebb) coefficient 
of the result is $c$. Note that we can choose this coefficient freely, as
the zeroth order \cheb and Taylor polynomials are constant.

The problem as stated has 14 boundary conditions, however
 because of conservation of energy, we only need 13 constraints. 
We choose the parameters as follows. Set
\begin{equation}\label{constrOne}
(d_1,d_2,d_3,d_4)=T_{SR}(x_2(-1),p_2(-1),y_2(-1),q_2(-1))-(x_0(1),p_0(1),y_0(1),q_0(1))\,,
\end{equation}
so that four boundary conditions become $d_1=d_2=d_3=d_4=0$.
Then we choose
\begin{itemize}
\item $c_1=c_2=c_5=c_6=0$, which imply that $x_1(0)=0\,, y_1(0)=0,x_2(0)=0\,, y_2(0)=0$.
\item $c_3,c_4,a_0,b_0$ in order to satisfy the equality
$$(\hat x_0(-1),\hat p_0(-1),\hat y_0(-1),\hat q_0(-1))=T_{PR}(x_1(1),p_1(1),y_1(1),q_1(1))\,.$$
\item $a_1=p_1(0)+d_2$
\item $a_2=p_2(0)+d_3$
\item $b_1=\pm\sqrt{8(1-\mu)-p_1^2(0)}$
\item $b_2=\pm\sqrt{8\mu-p_2^2(0)}$. 
\end{itemize}
Observe that we have to know in advance whether $q_1(0)$ and $q_2(0)$
are positive or negative, in order to choose the correct sign for $b_1$ and $b_2$.  
These choices are determined from the numerical approximation
of the solution we are trying to validate.  

We are ready to define the maps
$\XX=\AA^{8}\times\BB^4\times\mathbb{R}$ and  $\FF:\XX\to\XX$ by
\begin{equation}\label{eq:FF}
\FF(x_1,y_1,x_0,y_0,x_2,y_2,p_1,q_1,p_0,q_0,p_2,q_2,T_0)=(\tilde x_1,\tilde y_1,\tilde x_0,\tilde y_0,\tilde x_2,\tilde y_2,\tilde p_1,\tilde q_1,\tilde p_0,\tilde q_0,\tilde p_2,\tilde q_2,\tilde T_0)
\end{equation}
where
\begin{equation}\label{eq:fpe}
\begin{cases}
(\tilde x_1,\tilde y_1,\tilde x_0,\tilde y_0,\tilde x_2,\tilde y_2)= D^{-1}_{c_i}(T_1p_1,T_1q_1,T_0p_0,T_0q_0,T_2p_2,T_2q_2)\\
\tilde  p_1=T_1D^{-1}_{a_1}f_1(x_1,y_1,p_1,q_1)\\
\tilde  q_1=T_1D^{-1}_{b_1}g_1(x_1,y_1,p_1,q_1)\\
\tilde  p_0=T_0D^{-1}_{a_0}f_0(x_0,y_0,p_0,q_0)\\
\tilde  q_0=T_0D^{-1}_{b_0}g_0(x_0,y_0,p_0,q_0)\\
\tilde  p_2=T_2D^{-1}_{a_2}f_2(x_2,y_2,p_2,q_2)\\
\tilde  q_2=T_2D^{-1}_{b_2}g_2(x_2,y_2,p_2,q_2)\\
\tilde T_0=T_0+d_4\,,
\end{cases}
\end{equation}
and the constants $\{a_i,b_i,c_i,d_i\}$ are chosen as described above. 
It is straightforward to check that a fixed point $X$ of $\FF$ corresponds 
to an ejection-collision solution; in particular $d_2=d_3=d_4=0$, and because 
of the conservation of energy $d_1=0$ as well.
Note also that a solution of this equation is continuous and 
a-priori piecewise analytic, thanks to the 
decay rates imposed by the norms on the sequence spaces.  Moreover, since the 
CRTBP vector field of Equation \eqref{eq:PCRTBP} is real analytic on its domain of 
definition, the entire ejection-collision solution curve -- when projected 
back into rotating coordinates -- is globally analytic
up to the singularities at the endpoints.

\section{Fixed point problem for studying the bifurcations}\label{sec:bif}
We begin by discussing the local analysis near the bifurcation points
appearing in Theorems \ref{th:fixedC},\ref{th:muonehalf}, and \ref{th:muearthmoon}.
The argument is based on a procedure introduced in \cite{AGK}, which we now briefly review.
Note that the map $\mathcal{F}$ defined in Equation \eqref{eq:FF}
depends on the parameters $\mu$ and $C$, through 
the component functions $f_{0,1,2}$, $g_{0,1,2}$ of the vector fields 
defined in Equations \eqref{eq:PCRTBP}, \eqref{eq:regularizedSystem_m1}, 
and \eqref{eq:regularizedSystem_m2}.
We will vary one or the other of these parameters, 
which we refer to as $\tau$.  More precisely, $\tau = C$ in Theorems 
\ref{th:muonehalf}, and \ref{th:muearthmoon}, while $\tau = \mu$ 
in Theorem \ref{th:fixedC}.

We write $\FF_\tau$ to stress the explicit dependence 
of the map on the chosen parameter and, in a slight abuse of notation, 
define $\FF:\mathbb{R}\times\XX\to\XX$ by
\begin{equation}
\FF(\tau,X)=\FF_\tau(X)-X.
\label{fixedPb}
\end{equation}
Then, for a fixed value of $\tau$, $X_\tau$ is an ejection-collision orbit
with system parameter $\tau$ if and only if $\FF(\tau, X_\tau)=0$.
We write either $X(\tau)$ or $X_\tau$ depending on whether
 we wish to stress that $X$ is a function of 
$\tau$, or a branch of solutions parameterized by $\tau$.

Our goal is to numerically approximate the parameter value $\tau_*$ where 
a bifurcation occurs, and then to validate that there is 
a true bifurcation nearby via a contraction mapping argument.  
We also need mathematically rigorous bounds on the local branch(es) 
near the bifurcation point.  Our strategy is to develop a kind of 
normal form $g(\tau, \lambda)$ for the fixed point operator $\FF_\tau$, in such a way that the 
zero level set of $g$ describes the branch in a neighborhood
of the numerically computed bifurcation point.

To formalize the discussion, we begin
with a two-parameter family
of functions $X(\tau,\lambda)$ solving
\begin{equation}
(\id-\ell)\FF\bigl(\tau,X(\tau,\lambda)\bigr)=0\,,\qquad
\ell X(\tau,\lambda)=\lambda\hat X,
\label{ualphalambda}
\end{equation}
where $\id$ is the identity map on $\XX$, and 
$\ell:\XX\to\XX$ is a one-dimensional 
projection whose image approximates the kernel of $D \FF(\tau, X)$, 
 and $\hat X\in\XX$ is in the range of $\ell$.
Choose $\{e_j\}_{j\in\mathbb{N}}$ a basis of $\XX$, so that 
each of the functions $e_j$ is in either
$\AA_\rho$ or $\BB_\varrho$ depending on context.  Take
$\hat X=\sum_{j=0}^N \hat X_{j}e_j$ an approximate 
eigenvector of $D\FF(\tau,\bdot)$ corresponding to the numerically smallest eigenvalue.
We fix the scale by choosing $\hat X$ with $\sum_{j=0}^N \hat X_{j}^2=1$, and 
define the projection $\ell$ by
\begin{equation}
\ell X=\ell_0(X)\hat X\,,\qquad \ell_0(u)=\sum_{j=0}^N X_{j}\hat X_{j}.
\end{equation}
Here $X_{j}$ are the coefficients of $X$ in the basis $\{e_j\}$.

The goal now is to show that equation (\ref{ualphalambda}) 
has a smooth and locally unique solution
$X: I\times J\to\XX\,$ on some rectangle $I\times J$ in the parameter space.
Toward this end, define the bifurcation
function $g \colon I \times J \to \mathbb{R}$ by 
\begin{equation}
g(\tau,\lambda)=\ell_0\FF\bigl(\tau,X(\tau,\lambda)\bigr).
\label{gDef}
\end{equation}
We observe that $g(\tau, \lambda)$ is zero when $X(\tau, \lambda)$ is
a solution of $\FF(\tau,X)=0$.

We expand the coefficients of $X$ (and hence $X$ itself)
about the point $(\tau_0, \lambda_0)$ 
as Taylor polynomials in $\tau,\lambda$ of the form 
\begin{equation}
X(\tau,\lambda)=\sum_{j} X_j(\tau, \lambda) e_j\,,\quad
X_j=\sum_{0\le l+m\le M} X_{jlm}\left(\tau-\tau_0\over \tau_1\right)^l
\left(\lambda-\lambda_0\over\lambda_1\right)^m.
\label{taylortwo}
\end{equation}
Here $\tau_0,\tau_1,\lambda_0\lambda_1,M$ are as given in Table 1
(note that $\tau=\mu$ in the case $i=1$ and $\tau=C$ in the cases $i=2,3$).
\begin{table}[htp]\label{tab:bif}
\begin{center}
\begin{tabular}{|c|c|c|c|c|c|}
\hline
$i$ & $\tau_0$ & $\tau_1$ & $\lambda_0$ & $\lambda_1$ & $M$\\
\hline
1 & $1631\cdot 2^{-13}$ & $5\cdot 2^{-14}$ & $5\cdot 2^{-7}$ & $2^{-4}$&14\\
2 & $8453\cdot 2^{-12}$ & $ 2^{-12}$ & $51\cdot 2^{-8}$ & $25\cdot2^{-9}$&12\\
3 & $8397325\cdot 2^{-22}$ & $2^{-26}$ & $-2102\cdot 2^{-13}$ & $3\cdot2^{-13}$&8\\
\hline
\end{tabular}
\vskip2mm
\caption{Parameters for the bifurcations}
\end{center}
\end{table}
We refer to such expansions as $\XX-$Taylor or (XT) 
series in $\tau$ and $\lambda$, and to the truncations as XT polynomials.
Solving Equation (\ref{ualphalambda}) for the unknown branch
$X=X(\tau,\lambda)$ is equivalent
to finding a fixed point of the map $\FF_{\tau,\lambda}\,$, defined by
\begin{equation}
\FF_{\tau,\lambda}(X)=(\id-\ell)\FF_\tau(X)+\lambda\hat X.
\label{Fbetalambda}
\end{equation}
Our goal is to use
the Banach Fixed Point Theorem to obtain a true solution near our numerical approximation.

The map $\FF$ is not a (local) contraction, but it is compact.
More precisely, its derivative $D\FF$ has all eigenvalues, but finitely many,
of modulus less than 1. To prove the existence of a fixed point of $\FF$
by means of the Banach Fixed Point Theorem
we define the new mapping
\begin{equation}
\MM_{\tau,\lambda}(h)=\FF_{\tau,\lambda}(\bar X+ \Lambda_{\tau,\lambda} h)-\bar X+M_{\tau,\lambda} h\,,\qquad\Lambda_{\tau,\lambda}=\id-M_{T,\lambda},
\label{contractionTwo}
\end{equation}
where $\bar X$ is approximately fixed by $\FF_{\tau_0,\lambda_0}$, 
and $M_{\tau,\lambda}$ is a finite rank operator
having that $\Lambda_{\tau,\lambda}=\Id-M_{\tau,\lambda}$  
approximately inverts the operator $\Id-D\FF_{\tau_0,\lambda_0}(\bar X)$,
and we observe that, if $h$ is a fixed point for $\MM_{\tau,\lambda}$,
then $\bar X+ \Lambda_{\tau,\lambda} h$ is a fixed point for $\FF_{\tau,\lambda}$.
One can see that any solution, by the Implicit Function Theorem, 
depends analytically on $\tau$ and $\lambda$.

Let $D_r(z)$ denote the complex disk of radius $r$ centered at $z$.
We can, with the aid of a computer, prove the following theorem.  
See  Section \ref{sec:cap} for discussion of the implementation.
\begin{lemma}\label{capTwo}
For $i=1,2,3$, let $\tau_0,\tau_1,\lambda_0,\lambda_1,M$ 
as in Table 1, define $I=D_{\tau_1}(\tau_0)$ and $J=D_{\lambda_1}(\lambda_0)$.
There exists an XT polynomial $\bar X(\tau,\lambda)$  as in \eqref{taylortwo}, and
positive constants $\eps,r,K$ satisfying $\eps+Kr<r$, so that
\begin{equation}\label{bifcap}
\|\MM_{\tau,\lambda}(0)\|\le\eps\,,
\qquad\|D\MM_{\tau,\lambda}(v)\|\le K.
\end{equation}
for all $v\in B_{r}(0)$  and all
$\tau\in I$, $\lambda\in J$.
\end{lemma}

Combining the Banach Fixed Point Theorem with Lemma \ref{capTwo}, we now have
the following result:

\begin{proposition}\label{biflet}
For $i=1,2,3$, let $\tau_0,\tau_1,\lambda_0,\lambda_1,M$ as in Table 1, 
let $I=D_{\tau_1}(\tau_0)$ and $J=D_{\lambda_1}(\lambda_0)$.
Then for each $(\tau,\lambda)$ in $I\times J$, 
Equation \eqref{ualphalambda} has a unique solution $X(\tau,\lambda)$ in
$B_r(\bar X(\tau,\lambda))$.  Moreover, the map
$(\tau,\lambda)\mapsto X=X(\tau,\lambda)$ is analytic and, 
for each real $\tau\in I$, a function $X$ in $B\cap\ell^{-1}(J\hat X)$
is a fixed point of $F_\tau$ if and only if $X=X(\tau,\lambda)$
for some real $\lambda\in J$, with $g(\tau,\lambda)=\ell_0 X(\tau,\lambda)=0$.
\end{proposition}

\begin{figure}[ht]
\begin{center}
{\includegraphics[height=45mm, width=44mm]{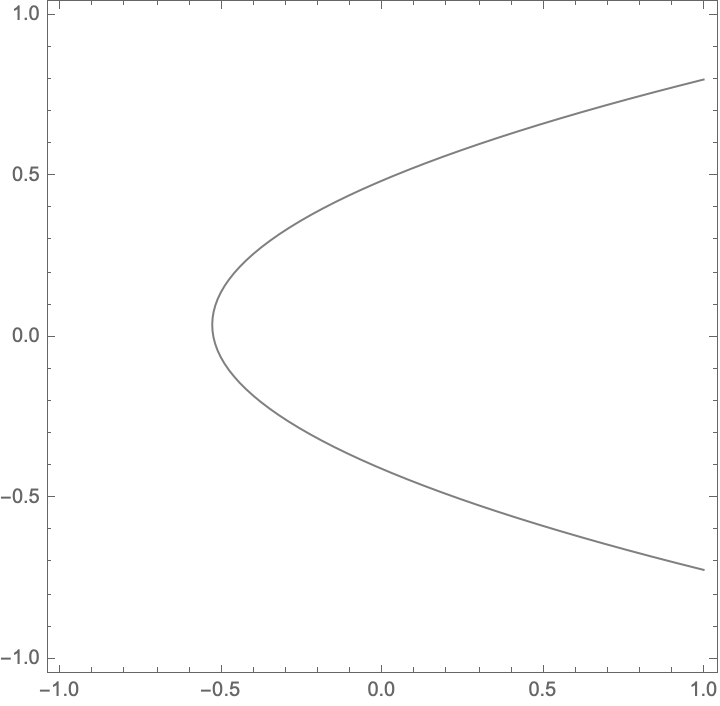}}
{\includegraphics[height=45mm, width=44mm]{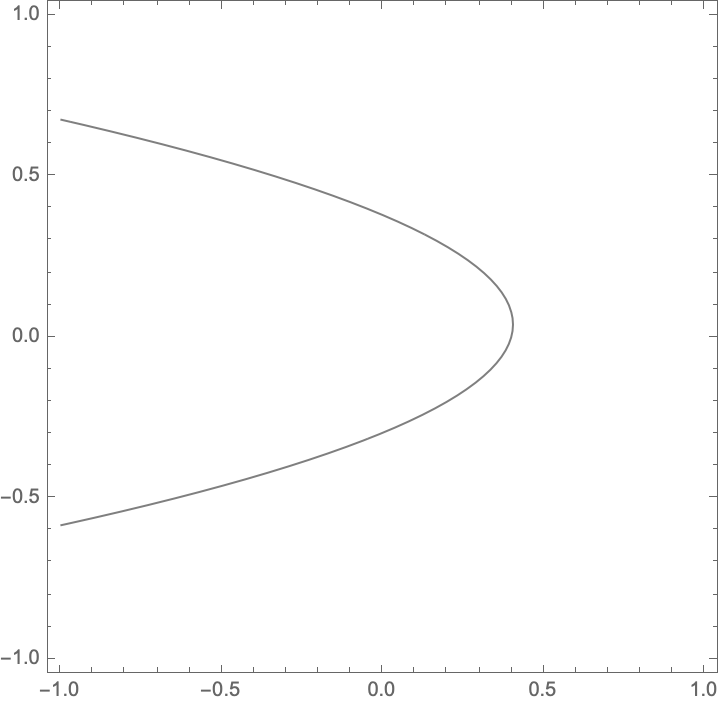}}
{\includegraphics[height=45mm, width=44mm]{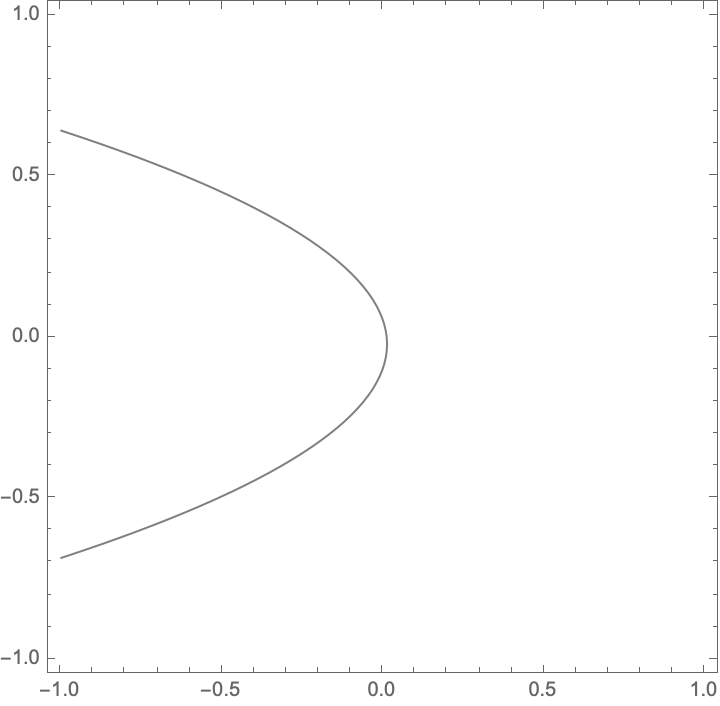}}
\caption{Approximate zero level curve for
 the functions $g(\tau,\lambda)$ describing the saddle-node bifurcations
in Theorems \ref{th:fixedC},\ref{th:muonehalf}, and \ref{th:muearthmoon}.}
\end{center}
\label{fig:bifsG}
\end{figure}

Figure \ref{fig:bifsG} depicts the
numerically computed zero level sets of the functions $g$ 
obtained in Proposition \ref{biflet}.
Based on the images depicted in Figure \ref{fig:bifsG}, one can guess that there are
saddle-node bifurcations in each of the three cases being considered.
We are left with the problem of understanding the bifurcation points in terms 
of the behavior of the zero level set of $g$.
Following  \cite{AK1}, a set of conditions sufficient for establishing the existence
of a saddle-node bifurcation are given below. 
The interested reader is referred to \cite{AK1} for the proof.

\begin{lemma}\label{lem:fold} {\rm (saddle-node bifurcation)}
Let $I=[\tau_1,\tau_2]$ and $J=[b_1,b_2]$.
Suppose that $g$ is a real-valued $C^3$ function
on the open neighborhood of $I\times J$, and let
$\dot g$ and $g'$ denote partial dedifferentiation 
with respect to the first and second argument, respectively.

Assume that 
\halign{\ #&\ $#$\hss&\quad#&\ $#$\hss&\quad#&\ $#$\hss\cr
(1) & g''>0 ~~on~~ I\times J, &
(2) & \dot g<0 ~~on~~ I\times J, & & \cr
(3) & g(\tau_1,0)\pm\half b g'(\tau_1,0)>0, &
(4) & g(\tau_2,\pm b)>0, &
(5) & g(\tau_2,0)<0. \cr}
\noindent
Then the solution set of $g(\tau,\lambda)=0$ in $I\times J$
is the graph of a $C^2$ function $\tau=a(\lambda)$,
defined on a proper subinterval $J_0$ of $J$.
This function takes the value $\tau_2$ at the endpoints of $J_0\,$,
and satisfies $\tau_1<a(z_2)<\tau_2$
at all interior points of $J_0\,$, which includes the origin.
\end{lemma}

Finally, we establish the following lemma with computer assistance.

\begin{lemma}\label{capThree}
For $i=1,2,3$, consider the solution $X(\tau,\lambda)$ obtained in Proposition \ref{biflet}.
For any $X(\tau,\lambda)\in B_r(\bar X(\tau,\lambda))$ the function
$g(\tau,\lambda)=\ell_0X(\tau,\lambda)$ if $i=1$ or $g(\tau,\lambda)=\ell_0X(\tau,-\lambda)$ if $i=2,3$ satisfies the
assumptions of Lemma \ref{lem:fold}.
\end{lemma}

\section{Gluing of Branches}\label{sec:branches}

The next step is to follow the local branches established in the previous section
away from the bifurcation point.  For example,
to prove Theorem \ref{th:fixedC} it is necessary to establish that the 
solution branches $B_1$ and $B_2$ emerging from the saddle-node bifurcation 
at $\mu=0.19934\ldots$ continue smoothly all the way to the mass ratio $\mu=1/2$. 
More precisely, since Lemma \ref{capThree} establishes that the branches emerging from 
the saddle-node bifurcation extend up to $\mu=1633 \cdot 2^{-13}$, we want to 
show that the branches exist over the parameter interval $\mu \in \II=[1633\cdot2^{-13},1/2]$.
This is done by numerically computing a high order Taylor expansion of
$K$ solution arcs on $K$
sub-intervals of $\II$, proving the existence of a true solution near each of
the numerical approximations using a fixed point argument, 
and showing that the sub-intervals ``glue together'' in the proper way.

So, take $\{e_j\}$ again to be the basis introduced in 
the previous section, and expand the solution coefficients on each 
sub-interval as $\XX$ Taylor polynomials in $\mu$, given by 
\begin{equation}
X_m(\mu)=\sum_{j} X_j^m(\mu) e_j\,,\quad
X_j^m =\sum_{0\le l\le L_m} X_{jl}^m \left(\mu-\mu_m \over \tau_m\right)^l\,,
\label{taylorone}
\end{equation}
where $\mu_m, \tau_m$, and $L_m$ 
for $1 \leq m \leq K$ are discussed further below.
We numerically compute  XT
polynomials $\bar X_\mu^m$ for $1 \leq m \leq K$
having that $\FF_\mu(\bar X_\mu^m) \simeq \bar X_\mu^m$ for each $m$, and 
finite rank operators $M_\mu^m:\XX\to\XX$ such that $\id-M_\mu^m$ is an approximate 
inverse of $\id-D\FF_\mu(\bar X_\mu^m)$ for each $m$.

As in the previous section, to use the Banach Fixed Point Theorem
we introduce new maps $\MM_\mu^m$ as
\begin{equation}
\MM_\mu^m(h)=\FF_\mu(\bar X^m+ \Lambda_\mu h)-\bar X^m+M_\mu^m h\,,
\qquad\Lambda_\mu^m=\id-M_\mu^m, \quad \quad \quad 1 \leq m \leq K.
\label{contraction}
\end{equation}
For $r>0$ and $w\in\XX$, let $B_r(w)=\{v\in\XX\,:\,\|v-w\|<r\}$ denote the ball of 
radius $r$ about $w$.  The branch $B_1$ (resp $B_2$) is partitioned into K = 59 
(resp.  K= 63) subintervals
of various widths, chosen adaptively depending on how far we are from the bifurcation points.
The exact values of the centers $\mu_m$ and widths $\tau_m$ of the subintervals, 
as well as the degrees of the Taylor expansions and the numerically computed
coefficients, are available in the data
files {\tt params.ads}, {\tt run\_b1.adb}, {\tt run\_b2.adb}.
The following lemma is then established with computer assistance.  See also Section \ref{sec:cap}.

\begin{lemma}\label{capOne}
The following holds for each value of $\mu_m,\tau_m,L_m$ in the data files mentioned above.
For each $1 \leq m \leq K$
there exist a XT polynomial $\bar X^m(\mu)$ of degree $L_m$ as described in \eqref{taylorone},
a bounded linear operator $M_\mu^m$ on $\XX$,
and positive real numbers $\eps_m, r_m, K_m$ satisfying $\eps_m+K_mr_m<r_m$, such that
\begin{equation}
\|\MM_\mu^m(0)\|\le\eps_m\,,\quad
\|D\MM_\mu^m(Y)\|\le K_m\,,\quad \forall Y\in B_{r_m}(0)
\label{apprfp}
\end{equation}
holds for all $\{\mu\in\complex: |\mu-\mu_m|<\tau_m\}$. 
\end{lemma}

So, in the case of the proof Theorem \ref{th:fixedC},
Lemma \ref{capOne} establishes the existence of 
122 subbranches of ejection collision orbits.
The next step is to check that the subbranches coincide at their endpoints, 
 all the way to the saddle-node bifurcation. This is proven with computer assistance, 
by showing that the solutions match at the end points of the sub-branches.  Then
they are, in fact, pairwise on the same branch. 

So, for each XT representation $X_\mu^m$ of an arc on the branch, let $P(X_\mu^m)$ 
be an enclosure of the image of the $(x_2,y_2,p_2,q_2)$ components of $X_\mu^m$, 
evaluated at $T=-1$.
\begin{lemma}\label{patchOne}
Suppose that a pair of adjacent solution arcs $X^{m}_\mu,X^{m+1}_\mu$ intersect at some 
value of $\bar\mu$, and assume that we have either
$$P(X^{m}_{\mu})\subset P(X^{m+1}_{\mu})\text{ or }P(X^{m+1}_{\mu})\subset P(X^m_{\mu}).$$
Then $X^m_{\mu}$ and $X^{m+1}_{\mu}$ represent different arcs on the same
branch of ejection-collision solutions.
\end{lemma}

\noindent Repeated use of Lemma \ref{patchOne} on each of the $K$ subintervals completes the proof.

We still need to prove that the branches $B_1$ and $B_2$ are connected to the
small branches emerging from the saddle-node bifurcation proved in Lemma 
\ref{capThree}.
To do so, let $\lambda_1=3/8+B_{2^{-8}}(0)$ and $\lambda_2=-343/512+B_{2^{-8}}(0)$, 
with $B_r(0)=[-r,r]$.
Let $X_{sn}$ be the solution obtained in Proposition \ref{biflet} 
when $i=1$, evaluated at $\mu=\bar\mu=1633\cdot2^{-13}$, 
and let $P_j(X_{sn})$, $j=1,2$ be an enclosure of the value of the 
$(x_2,y_2,p_2,q_2)$ components of $X_{sn}$, evaluated at $T=-1$, 
and $\lambda=\lambda_j$. Let $P_j(X_b)$, $j=1,2$, denote an enclosure of the 
values of the $(x_2,y_2,p_2,q_2)$ components of the two solutions considered in 
Lemma \ref{capOne} corresponding to $\mu=\bar\mu$, evaluated at $T=-1$. Then the values $P_j(X_{sn})$ are an enclosure of a point of the solutions belonging to the branches connected to the saddle-node bifurcation, while the values $P_j(X_b)$ are an enclosure of a point of the solutions belonging to the larges branches.
\begin{lemma}\label{patchTwo}
For $i=1,2$
$$P_j(X_b)\subset P_j(X_{sn})\quad\text{ and }P_1(X_b)\cap P_2(X_b)=\emptyset\,.$$
\end{lemma}
This implies that the two solutions obtained in Lemma \ref{patchOne} at $\mu=\bar\mu$ are different, and coincide with the solutions obtained in Proposition \ref{biflet} and Lemmas \ref{capThree}.
These results, together with the Contraction Mapping Theorem and the Implicit Function Theorem,
imply Theorem \ref{th:fixedC}.

\section{Computer representation and estimates}\label{sec:cap}

The computer assisted arguments discussed in Sections 
\ref{sec:bif} and \ref{sec:branches} are implemented in the
programming language Ada \cite{Ada}.  The remainder of the 
present section will serve as a rough guide for the reader 
interested in running these programs.  
The programs themselves are found at \cite{filesR3B}.

In the programs, many of the most important bounds are
obtained by enclosing the ranges of functions.  
More precisely, the function  $f:\XX\to\YY$
is inclosed by a function $F$ if assigning $F$ to each 
set $X\subset\XX$ of a given type ({\tt Xtype}) a set $F(X) = Y\subset\YY$
of a given type ({\tt Ytype}), in such a way that
$y=f(x)$ belongs to $Y$ for all $x\in X$.
In Ada, such a bound $F$ can be implemented by defining
a {\tt procedure F(X\,:\, in Xtype\,;\, Y\,:\, out Ytype)}.

We say that a number $x \in \mathbb{R}$ is representable if it can be 
expressed exactly as a floating point number in the computer.  
In our programs, 
a ball in a real Banach algebra $\XX$ with unit ${\bf 1}$
is represented using a pair {\tt S=(S.C,S.R)}.  Here
{\tt S.C} denotes a representable number ({\tt Rep})
while {\tt S.R} represents a nonnegative representable number ({\tt Radius}).
Then, for example, a ball in $\XX$ is represented by the set
 $\langle{\tt S},\XX\rangle=\{x\in\XX:\|x-({\tt S.C}){\bf 1}\|\le{\tt S.R}\}$.

In the case that $\XX=\mathbb{R}$, we refer to the data type just described by {\tt Ball}.
Numerical implementation of the enclosures of standard functions on the data 
type {\tt Ball} are defined in the packages {\tt Flts\_Std\_Balls}.
The packages {\tt Vectors} and {\tt Matrices} extend the basic functions, 
facilitating the enclosure of vector and matrix valued functions.  
These kinds of data structures are used frequently in 
computer-assisted proofs, so we will describe only the
more problem-specific features of our implementation.

In the present work, away from bifurcation points, 
 one parameter families (i.e. branches) 
of solutions are described by one variable Taylor series.
Two variable power series are used near bifurcation points.
The one variable and two variable Taylor series are represented
using data types called {\tt Taylor1} and  {\tt Taylor2} respectively.
These types consists of one and two variable power series, whose 
coefficients are of type {\tt Ball}.
Definitions, basic properties, and procedures are found in the packages {\tt Taylors1}
and  {\tt Taylors2} respectively.
The interested reader will find these packages at \cite{brus}.

Similarly, functions in $\AA$  are represented using a type
called {\tt Cheb}, defined in a package called {\tt Chebs}.
This data type has coefficients in some Banach algebra with unit $\XX$.
For example, in the present work our {\tt Cheb} objects have
coefficients of either type {\tt Taylor1} or {\tt Taylor2}, 
depending on the application.  
The {\tt Cheb} data type is represented by a triple {\tt F=(F.C,F.E,F.R)},
with {\tt F.C}  an {\tt array(0..K) of Ball},
{\tt F.E}  an {\tt array(0..2*K) of Radius}, and {\tt F.R}
an object of type {\tt Radius}.
The later describes the ellipse of analyticity 
on which the functions in $\AA$ are defined.  That is {\tt F.R}$\,=5/4$.
The set $\langle{\tt F},\AA\rangle$
consists of function $u=p+h\in\AA$ with
$$
p(t)=\sum_{j=0}^K \langle{\tt F.C(J)},\XX\rangle\,T_j(t)\,,\qquad h=\sum_{j=0}^{2K}h^{j}\,,\qquad
h^{j}(t)=\sum_{m\ge j} h^{j}_{m}\,T_m(t).
$$
Here $\|h^{J}\|\le{\tt F.E(J)}$, for all $J$.

The operations needed in our proof  
are bound efficiently using these enclosures.
The data types {\tt Chebs} and {\tt Taylors1} are
used to define the type {\tt XX}, which is in turn used
to represent the subspaces $\proj_{\{j\}}\AA_\BB$ (coefficient modes)
and $\proj_{\{j,j+1,\ldots\}}\AA_\BB$ (error modes).
Partitions of unity (in the sense of direct sums) are defined using 
Arrays of {\tt CylinderMode}.   The procedure {\tt Make}
implements these partitions 
for the spaces $\AA_j$ described in Section 6.

Operations needed in the implementation of a Newton-like 
operator are defined in the higher level packages 
{\tt Linear} and {\tt Linear.Contr}, which have been used to solve a
wide variety of problems.  
The procedures defined in these packages are designed to work 
with generic function types {\tt Fun}, as long as their associated {\tt Modes}
are specified.  In the present work, these packages are used to implement
quasi-Newton maps $\NN$,  estimate derivatives (such as the
linear operator such as $D\NN(h)$) and to check contraction bounds like those
in \eqref{bifcap} and \eqref{apprfp}.
A much more complete description of of the package {\tt Modes}
is found for example in \cite{MR3886634,MR4075873}.

\appendix

\section{Remarks on the literature} \label{sec:literature}

The first computer assisted proofs dealing with the 
dynamics of the planar CRTBP appeared in the papers
\cite{MR1947690, MR2112702} 
by the first author, and established the existence of
branches of periodic orbits, as well as chaotic 
subsystems with positive topological entropy for the 
CRTBP with equal masses.  
Shortly thereafter, Wilczak and Zgliczy\'{n}ski 
generalized these techniques to prove the existence of 
heteroclinic connections between periodic orbits in the Sun-Jupiter system
\cite{MR1961956, MR2174417}.

Over the last 20 years, many additional
theorems have been proven for planar and spatial CRTBPs
using computer assisted arguments.
Rather than attempting a systematic review, we 
refer the interested reader to 
Section 1.3 of \cite{MR4576879}, where a more thorough discussion
and many additional references are found.  We also mention the general 
review article of van den Berg and Lessard 
\cite{jpjbReview} on computer assisted methods of proof in 
dynamical systems theory, and the review of G\'{o}mez-Serrano 
\cite{MR3990999} on computer assisted proofs for PDEs.
Undergraduate and graduate level discussions of validated numerical techniques,
appropriate for use in computer assisted proofs, are found in the
book of Tucker \cite{MR2807595}, in the review article of Rump 
\cite{MR2652784}, 
and also in the monograph of Nakao, Plum and Watanabe \cite{MR3971222}.

A central question in nonlinear analysis is to understand 
the dependence of solutions on problem parameters.  Given the importance of the question,
it is natural that a number of authors have devoted substantial effort 
to the development of computer assisted methods of proof for
 continuation and bifurcation arguments. 
Perhaps the first computer assisted result of this kind was the 
paper of Plum \cite{MR1354655} on continuation of branches 
of solutions of one parameter families of second order scalar BVPs.  
Computer assisted proofs for bifurcations in the Rayleigh-B\'{e}nard
problem, using similar techniques, 
 were given by Nakao, Watanabe, Yamamoto, Nishida, and Kim
in \cite{MR2639642}.

The paper \cite{AK1}, by Koch and the first author,
provides a coherent example of how techniques of 
validated continuation can be combined with
mathematically rigorous bifurcation  
analysis to obtain mathematically rigorous bifurcation
diagrams.  The authors studied
 equilibrium solutions for a one parameter family of
 scalar PDEs defined on a one dimensional spatial 
 domain with Dirichlet boundary conditions.  An example 
 application to a scalar PDE defined on a
 two dimensional spatial domain with periodic boundary conditions was given 
by van den Berg and Williams in \cite{MR3636312}.
For other extensions and applications of validated continuation 
we refer to the papers
of Day, Lessard, and Mischaikow \cite{MR2338393}, 
Gameiro, Lessard, and Mischaikow \cite{MR2487806},
and to the references therein.  

More recently, Gazzola, Koch and the first author developed a computer
assisted continuation methods based on a 
single variable Taylor expansion for
branches of  solutions of nonlinear problems.  See \cite{AGK} for an application
to the planar Navier-Stokes equation. 
The authors set up the continuation problem as a 
fixed point problem on a bigger Banach space, and solve ``all-at-once'' for a high order
expansion of the branch using any desired basis functions.
A similar all-at-once set up, employing a 
multi-variable Taylor expansion of a suitable normal form, is used to prove
the existence of pitchfork bifurcations.
The approach of the paper just cited
is the basis of the high order continuation and bifurcation
results of the present work. 

Breden and Kuehn developed a similar computer assisted 
continuation method using general expansions for branches of 
solutions of nonlinear problems \cite{MR4077213,MR4595838} (not just Taylor). 
The idea behind their work just cited is that, when studying random perturbations
of differential equations, it is important to expand the branch using
basis functions adapted to the distribution of the noisy parameter.  
We also mention the recent work 
of Calleja, Garcia-Azpeitia, H\'{e}not, Lessard, and the 
second author in \cite{calleja2024lagrangetrianglefigurechoreography},
where the existence of a global branch of choreographic 
solutions is established using computer assisted methods of proof.
The branch of periodic orbits are expressed using 
Fourier-Chebyshev series. More precisely, each choreography along the 
branch is represented by a Fourier series and the dependence of the 
branch the parameter is represented by a Chebyshev series (rather than 
multiple Taylor series as in the present work).
One quarter of the global branch of solutions is obtained using a single 
Chebyshev domain, and the remainder of the global branch 
is obtained by a symmetry argument.

We also mention the work  
of Walawska and Wilczak \cite{MR3923486}, 
where the authors use computer assisted
continuation and bifurcation arguments to study families of 
periodic orbits in the spatial CRTBP.  More than this, 
they follow a family of out-of-plane periodic orbits from a
symmetry breaking bifurcation, through a number of 
secondary bifurcations, and back to where it began.
That is, they prove the existence of a closed circle 
-- or global family -- of periodic orbits locally parameterized by energy 
(note that a circle of periodic orbits is topologically a torus).
Note also that in the work of \cite{MR4576879}, already mentioned 
above, the authors prove theorems on the existence of one parameter
families of periodic orbits -- parameterized by energy -- which
``pass through'' collision.

As a final remark we note that while the present work and 
\cite{MR4576879} both deal with computer assisted existence 
proofs for ejection-collision orbits,
the two works have very different agendas and there are 
important differences in their approaches. 
In particular, the results of \cite{MR4576879} focus on 
fixed values of the mass ratio and energy.  The results  
are based on the study of certain multiple shooting 
problems formulated in $\mathbb{R}^d$ (usually with $d$ large).
A-posteriori analysis of these shooting problems requires only the 
application of an elementary finite dimensional
Newton-Kantorovich argument. 
The techniques lead to transversality of the results,
leading to only local information about the existence of branches.  

The goal of that work was to develop ``plug and play'' shooting templates, suitable 
for a wide variety of computer assisted proofs involving any number of 
collision or near collision events.  Because of this, all
infinite dimensional complications are
hidden behind the sophisticated CAPD library for computing 
mathematically rigorous set enclosures of 
solutions to initial value problems and 
variational equations.  An excellent overview of the 
CAPD library, with many worked examples and 
a wealth of additional references is found in the review article of 
Kapela, Mrozek, Wilczak, and Zgliczy\'{n}ski \cite{CAPD_paper}.

Alternatively, the present work 
formulates a single fixed point equation describing an
ejection-collision orbit or a branch of ejection-collision orbits
 in an appropriate product of function spaces.
The advantage of the functional analytic approach in the
 present context is that the problem parameters appear 
 in the fixed point operator in a completely explicit way,
rather than entering the problem implicitly through 
flow map as in \cite{MR4576879}. 
This set-up facilitates high order
approximation of solution branches using techniques 
discussed in Section \ref{sec:bif}.
The proofs exploits software libraries developed by
Hans Koch and the first author over the course of the last 20 
years.  See for example the works of 
\cite{brus,MR3886634,MR4075873},
and also the references discussed therein as this list is 
by no means exhaustive.

\bibliographystyle{unsrt}
\bibliography{papers}

\end{document}